\documentclass[12pt]{article}
\usepackage{latexcad}
\usepackage[active]{srcltx}
\usepackage[centertags]{amsmath}
\usepackage{amsfonts}
\usepackage{amssymb}
\usepackage{amsthm}
\usepackage{graphicx}
\newcommand{\Author}[1]{{\bf \begin{center} {#1}
  \end{center}}\vspace*{1mm}}

\begin{document}
\begin{center}
{\Large \bf 4-critical wheel graphs of higher order.}
\Author{ Dainis ZEPS
\footnote{Author's address: Institute of Mathematics and Computer
Science, University of Latvia, 29 Rainis blvd., Riga, Latvia.
{dainis.zeps@lumii.lv}} }
\end{center}
\begin{abstract}

4-critical wheel graphs of higher order are considered concerning their belonging to free-planar or free-Hadwiger classes.
\end{abstract}

We are using terminology from \cite{moh 01,kra 94,ze 98}.

In \cite{jam1} Jose Antonio Martin Hernandez suggested to use the theory of free minor closed classes of graphs in trying to prove the Hadwiger conjecture. We present here some examination in this direction.

For $k>1$, let us call Hadwiger class of order $k$ $H_k=N_0(K_k)$, i.e. the one generated by the forbidden graph $K_k$. Then, following many researchers, Hadwiger conjecture may be formulated as follows: In the Hadwiger class of order $k>1$ graphs are $k-1$-colorable. For $k<7$ this assumption is proved to be true, but other cases remain  hypothetical.

In \cite{jam1} J.A.M. suggested to examine whether $k$-critical graphs belong to the class $Free(H_{k+1})$ \cite{kra 94,ze 98}. [Graph belonging to $Free(H_{k+1})$ we call free-Hadwiger graph.] Following Kratochvil theorem \cite{kra 94,ze 98}, $Free(H_{k+1})$ would be the class without minors $K_{k+1}^-$ and $K_{k+1}^\odot$, i.e. $K_{k+1}$ without edge or with split vertex.

It is easy to see that this assumption works for cases $k<4$.
Further, we present sequence of 4-critical graphs, $G_9$ (fig. 1), $G_7$ (fig. 2), $G_5$ (fig. 3), $G_3$ (fig. 4), where none of the graphs belong to the class $Free(H_5)$, i.e. they are not free-Hadwiger graphs because they all contain the minor $K_5^-$. Thus, the suggestion in \cite{jam1}, at least for case $k=5$,  could be declined. Nevertheless, we find it useful that help us to forward some new conjectures.
These newly discovered 4-critical graphs may be arranged in a sequence where each of them is a minor of the following, i.e.,
$$G_3 \prec G_5 \prec G_7 \prec G_9 .$$

Indeed, it is easy to see that $G_7$ may be received from $G_9$ by contracting two successive edges on the sides. The same applies to the inclusive pairs $G_5\prec G_7$ and $G_3\prec G_5$.

This sequence may be replenished with wheel graphs $W_i$, $i=3,5,7,9$, as minors of these graphs. See what we get in fig. 5. We use here almost trivial fact that $W_i\prec G_i$ for $i=3,5,7,9$.

All graphs $G_i$, $i>=3$ may be considered as first order higher wheels. All they have as minor $K_5^-$. Thus, they all have \emph{minor brackets}  $<K_5^-,K_5>$ and $<K_{3,3}^-,K_{3,3}>$ where minor bracket for graph $G$ we define as a pair of graphs $<a,b>$ where $a$ is minor of $G$, but $b$ isn't. On the other side, ordinary [zero order]wheels $W_i$, $i>4$ have minor brackets $<W_4,K_5^->$ and $<C_6^+,K_{3,3}^->$. Besides, all these minor brackets are simple, i.e., tightest possible ones in the very natural sense. Moreover, higher than first order wheels should have the same minor brackets as first order wheels.

Further, we raise the question 1) are the only 4-critical graphs that belong to $Free(H_5)$,  [zero order] wheel graphs and wheels with split edges [see fig. 6], i.e. 2) are all higher order wheels non-free Hadwiger graphs? Thus, does there exist a critical free-Hadwiger non-free-planar 4-critical graph?

Besides, we make some judgements concerning graph $G_3$. It is easy to see that $G_3$ is the 'cube with one corner cut off' graph. We could ask about higher order cube graphs, what we could get after cutting off some of cube's corners? For 4-cube, it is rather easy to see that cutting off corners [one, two, three] can't give any 5-critical graph. See, for example, one case of 4-cube graph with two corners of the cube cut off on fig. 7. Judging from the 4-cube graph experience, we conjecture that none of the higher order cube graphs with cut off corners of the cube can be a critical graph.

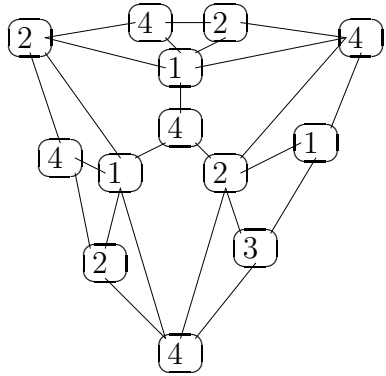
\placedrawing{HC5.LP}{Example of 4-critical graph that is not free-planar, not even free-Hadwiger graph. This is a first order higher wheel $G_9$.}{fig0}

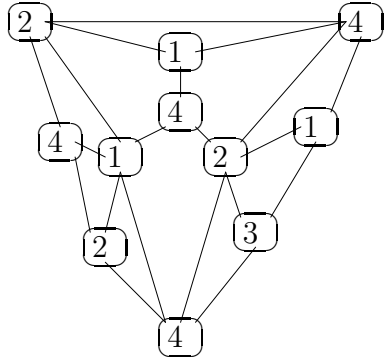
\placedrawing{HC4.LP}{First order higher wheel $G_7$. It is obtained from $G_9$ by contracting two of its side edges.}{fig1}

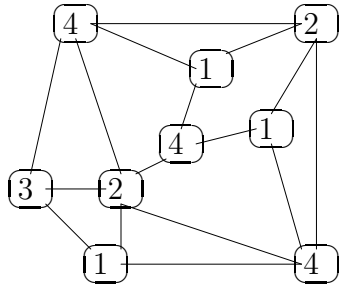
\placedrawing{HC3.LP}{First order higher wheel $G_5$. It is obtained from $G_7$ by contracting two of its side edges.}{fig2}

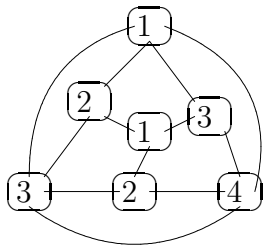
\placedrawing{HC2.LP}{First order minimal possible higher wheel $G_3$. It is obtained from $G_5$ by contracting two of its side edges. It is minimal in the sequence of these graphs obtained by contractions giving graph distinct from $K_4$. Next graph by contractions should be $K_4$. $G_3$ has minor $K_5^-$, thus, it is not Free-Hadwiger graph. Besides, it is easy to see that $G_3$ can be imagined as the cube with one corner cut off. }{fig3}

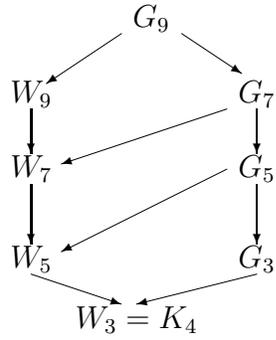
\placedrawing{HC1.LP}{Lattice of 4-critical graphs, where arrows show reductions of graphs by contractions of edges. In place of $G_9$ may stand any $G_i$ with odd $i$. Any next column in the lattice would be  wheel graphs with higher order.}{fig4}

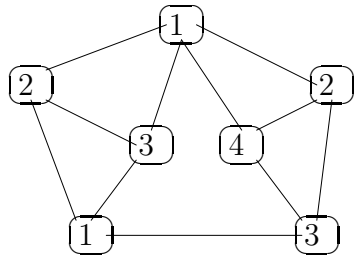
\placedrawing{HC6.LP}{4-critical graph obtained from $W_5$ by spitting of an edge [or spoke of the wheel]. The operation of the edge splitting preserving 4-criticality is easy to be generalized.}{fig5}

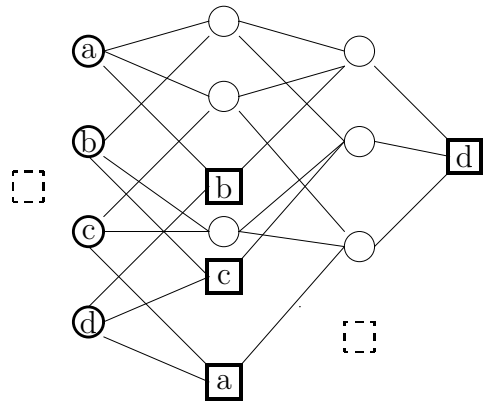
\placedrawing{HC7.LP}{4-cube with two corners [dash boxes] cut off. Marked vertices, circles and boxes correspondingly, are neighbors of cut off corners by hyperplane and thus all connected via edges [not drawn in the figure]. Letters denote distinct colors of vertices. It is easy to see that the graph is not 5-chromatic. }{fig6}


\begin{thebibliography}{99}
\bibitem{moh 01}
Mohar Bojan, Thomassen Carsten. {\em Graphs on Surfaces}, J. Hopkins Univ. Press, 2001.
\bibitem{jam1}
Jose Antonio Mart\'{\i}n Hernandez, {\em On Hadwiger conjecture}, personal e-mail communication, 2006.
\bibitem{kra 94}
Kratochv\'{\i}l J. {\em About minor closed classes and the generalization of the
notion of free-planar graphs}, personal communication, 1994, 2pp.
\bibitem{ze 98}
Zeps D. {\em Free Minor Closed Classes and the Kuratowski
Theorem}, KAM Series, 98-409, Prague, 1998, 10 pp.
\end{thebibliography}
\end{document}